\documentclass[10pt]{amsart}
\usepackage{mathptmx}
\usepackage{amsmath}     
\usepackage{amssymb}
\usepackage{array}
\usepackage{geometry}
\usepackage[bookmarks=true,colorlinks=true, pdfstartview=FitV, linkcolor=black, citecolor=blue, urlcolor=black]{hyperref}
\usepackage{movie15}

\usepackage{color}
\definecolor{DarkRed}{rgb}{0.55,.00,0.2}
\definecolor{DarkGrey}{rgb}{0.35,.35,0.35}

\theoremstyle{definition}

\theoremstyle{remark}

\numberwithin{equation}{section}

\begin{document}

\title{Discrete orthogonal polynomials associated with Macdonald function}

\author{{\bf S. Yakubovich} \\
 {\em {Department of Mathematics, Fac. Sciences of University of Porto,\\ Rua do Campo Alegre,  687; 4169-007 Porto (Portugal) }}}
  \thanks{ E-mail: syakubov@fc.up.pt} 
\vspace{3mm}

\thanks{ The work was  partially supported by CMUP, which is financed by national funds through FCT(Portugal),  under the project with reference UIDB/00144/2020. }

\subjclass[2000]{  33C10, 42C05, 44A15 }

\date{\today}


\keywords{Classical orthogonal polynomials, Charlier polynomials, Meixner polynomials, Laguerre polynomials, modified Bessel function}

\begin{abstract}  New sequences of discrete orthogonal polynomials associated with the modified Bessel function $K_\mu(z)$ or Macdonald function are considered.  The corresponding weight function is $\lambda^k \rho_{k+\nu+1}(t)/ k!$, where $\ k \in \mathbb{N}_0, \ t \ge 0,\ \nu > -1,\ 0 < \lambda < 1,\  \rho_{\mu}(z) = 2 z^{\mu/2} K_\mu\left( 2\sqrt z\right)$. The limit case $t=0$ corresponds to the Meixner polynomials. Various properties, differential-difference recurrence relations  are established.  The modified sequence of polynomials with the weight $\lambda^k \rho_{k+\nu+1}(\lambda t)/ k! $  is investigated as well.
\end{abstract}

\maketitle

\markboth{\rm \centerline{ S. Yakubovich}}{discrete orthogonal  polynomials}

\section{Introduction and preliminary results}

Let $\nu > -1,\   t \ge  0,\ 0 < \lambda < 1$ be parameters and consider the sequence of orthogonal polynomials $\{P_n\}_{n\ge 0}$ of degree $n$ on the integers $\mathbb{N}_0$, satisfying  orthogonality conditions

$$\sum_{k=0}^\infty P_n(k; t,\lambda) P_m(k; t,\lambda) \rho_{k+\nu+1}(t) { \lambda^k\over k!} =  \delta_{n,m},\quad 0 < \lambda < 1.\eqno(1.1)$$ 
Here $\delta_{n,m},\   n,m\in\mathbb{N}_{0}$ is the Kronecker symbol and $\rho_\mu(t)= 2 t^{\mu/2} K_\mu\left( 2\sqrt t\right)$, where $K_\mu(z)$ is the modified Bessel function or Macdonald function \cite{Bateman}, Vol. II.  It is represented by the integral \cite{Bateman}, Vol. II

$$   \rho_\mu(t) = \int_0^\infty  e^{-t/x -x}  x^{\mu -1} dx, \quad \mu > 0,\ t \ge 0.\eqno(1.2)$$
The choice of $\lambda$ in (1.1) is motivated by the convergence of the corresponding series since, evidently,  from (1.2) it has the inequality $  \rho_\mu(t)  \le \Gamma(\mu)$, where $\Gamma(z)$ is the Euler gamma function \cite{Bateman}, Vol. I.  Moreover,  as we see from (1.2) $\rho_{k+\nu+1}(0) = \Gamma(k+\nu+1)$, and the sequence $\{P_n (x; 0, \lambda) \}_{n\ge 0}$ in (1.1) is related to the Meixner polynomials \cite{Sze}.  

As an immediate consequence of the asymptotic behavior for  the modified Bessel function at infinity and near the origin as a function of $t$ (cf. \cite{Bateman}, Vol. II)  we find 
$$\rho_\mu (t)= O\left( t^{(\mu-|\mu|)/2}\right),\  t \to 0,\ \nu\neq 0, \quad  \rho_0(t)= O( \log t),\ t \to 0,\eqno(1.3)$$

$$ \rho_\mu(t)= O\left( t^{\mu/2- 1/4} e^{- 2\sqrt t} \right),\ t \to +\infty.\eqno(1.4)$$
This function  can be represented in terms of Laguerre polynomials (cf. \cite{PrudnikovMarichev}, Vol. II, Entry 2.19.4.13 )
$${(-1)^n t^n\over n!}\  \rho_\mu(t)=   \int_0^\infty x^{\mu+n -1} e^{-x - t/x}  L_n^\mu(x) dx,\quad    n \in\mathbb{N}_{0}.\eqno(1.5)$$
Further, it has a relationship with the  Riemann-Liouville fractional integral \cite{YaL}

$$ \left( I_{-}^\alpha  f \right) (t)  = {1\over \Gamma(\alpha)} \int_t^\infty (x-t)^{\alpha-1} f(x) dx,\quad  {\rm Re} \alpha > 0,\eqno(1.6)$$
namely,  we get  the  formula 
$$\rho_\mu(t)= \left( I_{-}^\mu \rho_0 \right) (t),\ \mu >0.\eqno(1.7)$$
Hence the index law for fractional integrals immediately implies

$$ \rho_{\nu+\mu} (t)= \left( I_{-}^\nu \rho_\mu \right) (t)=   \left( I_{-}^\mu \rho_\nu \right) (t).\eqno(1.8)$$
The corresponding definition of the fractional derivative presumes the relation $ D^\mu_{-}= - D  I_{-}^{1-\mu}$.   Hence for the ordinary $n$-th derivative of $\rho_\mu$ we find

$$D^n \rho_\mu(t)= (-1)^n \rho_{\mu-n} (t),\quad n \in \mathbb{N}_0.\eqno(1.9)$$
 Recalling (1.3) and integrating by parts, it is not difficult to establish  the following recurrence relation for $\rho_\mu$ 

$$\rho_{\mu+1} (t) =    \mu \rho_\mu(t)+ t \rho_{\mu-1} (t),\quad \mu \in \mathbb{R}.\eqno(1.10)$$
In the operator form it can be written as follows

$$\rho_{\mu+1} (t) = \left( \mu - tD \right) \rho_\mu(t).\eqno(1.11)$$

As is known,   up to a normalization factor  orthogonality conditions (1.1)  are equivalent  to the following equalities
$$\sum_{k=0}^\infty P_n(k; t,\lambda)  \rho_{k+\nu+1}(t) { \lambda^k k^j\over k!} = 0,\ j =0,\dots, n-1,\eqno(1.12)$$ 

$$\sum_{k=0}^\infty P_n(k; t,\lambda)  \rho_{k+\nu+1}(t) { \lambda^k (-k)_j  \over k!} = 0,\ j =0,\dots, n-1,\eqno(1.13)$$ 
where $(a)_j= a( a+1)\dots (a+j-1)$ is Pochhammer's symbol \cite{Bateman}, Vol. I.  As it follows from the theory of orthogonal polynomials \cite{Sze},  a sequence $\{p_n\}_{n\ge 0}$ satisfies the three term recurrence relation 

$$x p_n(x) = A_{n+1} p_{n+1}(x) +B_n p_n(x) + A_{n} p_{n-1}(x),\eqno(1.14)$$
where $p_{-1}(x) \equiv 0,\  p_n(x) = a_n x^n+ b_n  x^{n-1}+ \dots,$ and 

$$A_n= {a_{n-1} \over a_n},\quad B_n= {b_{n}\over a_n} - {b_{n+1}\over a_{n+1}}.\eqno(1.15)$$
As a consequence of (1.17) the Christoffel-Darboux formula takes place

$$\sum_{k=0}^n   p_k(x) p_k(y) = A_{n+1} \frac{  p_{n+1}(x) p_n(y) -  p_{n}(x) p_{n+1}(y)}{x-y}.\eqno(1.16)$$
Finally in this section, denoting by $\omega_k (t,\lambda) =  \rho_{k+\nu+1}(t) \lambda^k / k! $, we employ (1.10) to establish the Pearson-type equation  for the weight $\omega_k$

$$k(k-1) \omega_k   -  (k-1) (k+\nu) \lambda \omega_{k-1}  -   t \lambda^2  \omega_{k-2}  = 0.\eqno(1.17)$$
Differentiating by $t$ and recalling (1.9), we get

$${\partial \omega_k\over \partial t} = - {\lambda\over k} \  \omega_{k-1} ,\ k \in \mathbb{N}.\eqno(1.18)$$
A differentiation by $\lambda $ gives the equality 

$$ \lambda {\partial \omega_k\over \partial \lambda} = k \omega_k .\eqno(1.19)$$ 
Similarly,  via (1.9), (1.10) we derive  

$$ {\partial \over \partial \lambda} \left\{ \lambda^{\nu+1} \omega_k \right\} = (k+\nu+1) \rho_{k+\nu+1}(t) {\lambda^{k+\nu}\over k!}
= \rho_{k+\nu+2}(t) {\lambda^{k+\nu}\over k!} - \rho_{k+\nu}(t) { t \lambda^{k+\nu}\over k!}$$

$$= \lambda^{\nu-1} (k+1) \omega_{k+1}  -  {t \lambda^{\nu+1}\over k}  \omega_{k-1}, $$
i.e., invoking (1.18),  (1.19), we end up with the partial differential-difference  equation for the weight   $\omega_k$

$$  \lambda {\partial \omega_k \over \partial \lambda}  - t {\partial \omega_k\over \partial t} =  {\partial \omega_{k+1} \over \partial \lambda}  - (\nu+1) \omega_{k}.\eqno(1.20)$$
Differentiating through by $t$ in (1.20) and using (1.18), (1.19), we arrive at the ordinary second order  differential equation for the weight $\omega_k$

$$t  {d^2 \omega_k \over d t^2} - (k+\nu) {d \omega_k \over d t} - \omega_k = 0.\eqno(1.21)$$

\section{Orthogonal polynomials for the weight $\omega_k$}

In order to proceed our investigation we define the power moment $\mu_n(t)$ and factorial moment $\gamma_n(t)$  for the weight $\omega_k$ as follows

$$\mu_n = \sum_{k=0}^\infty  \rho_{k+\nu+1}(t) { \lambda^k k^n \over k!},\eqno(2.1)$$

$$\gamma_n = \sum_{k=0}^\infty  \rho_{k+\nu+1}(t) { \lambda^k \over (k-n)!}.\eqno(2.2)$$
Their values are given by the following lemma.

{\bf Lemma 1}. {\it Let $n \in \mathbb{N}_0,\ t > 0,\ 0 < \lambda < 1,\ \nu > -1$.  Moments $(2.1), (2.2)$ for the weight $\omega_k$ have the values

$$\mu_n  =   \sum_{j=0}^n S(n,j)\  {\lambda^j  \rho_{\nu+j+1} \left( t(1-\lambda) \right) \over (1-\lambda)^{\nu+j+1} },\eqno(2.3)$$

$$ \gamma_n  =  {\lambda^n   \rho_{\nu+n+1}\left( t(1-\lambda) \right) \over (1-\lambda)^{\nu+n+1}},\eqno(2.4)$$
where $S(n,j)$ are Stirling numbers of the second kind.}

\begin{proof}  In fact,  employing representation (1.2),  we have from (2.1)

$$\mu_n = \sum_{k=0}^\infty  \rho_{k+\nu+1}(t) { \lambda^k k^n \over k!} =   \int_0^\infty  e^{-t/x -x}  x^{\nu}  \sum_{k=0}^\infty  { (x\lambda)^k k^n \over k!} dx $$

$$=  \int_0^\infty  e^{-t/x -x}  x^{\nu}  \left( x {d\over dx} \right)^n \sum_{k=0}^\infty  { (x\lambda)^k  \over k!} dx =   \int_0^\infty  e^{-t/x -x}  x^{\nu}  \left( x {d\over dx} \right)^n \bigg\{ e^{\lambda x} \bigg\} dx,\eqno(2.5)$$
where interchange of the order of integration and summation and the multiple differentiation under the series sign in (2.5) are  permitted via the uniform convergence of the series with respect to $x \ge 0$ on closed intervals of $\mathbb{R}_+$ and the estimate

$$  \sum_{k=0}^\infty  \rho_{k+\nu+1}(t) { \lambda^k k^n \over k!} \le \sum_{k=0}^\infty  \Gamma(k+\nu+1) { \lambda^k k^n \over k!} < \infty,\quad  0 < \lambda < 1.$$
Meanwhile, denoting the differential operator by $D \equiv {d\over dx}$,  we appeal to the known operational identity (see \cite{Riordan}, Section 6.6)

$$ \left( x D \right)^n = \sum_{j=0}^n S(n,j) x^j D^j,\eqno(2.6)$$
involving  Stirling numbers of the second kind \cite{Riordan}.  Therefore, we find from (2.5) and (1.2) 

$$\mu_n =   \sum_{j=0}^n S(n,j) \lambda^j \int_0^\infty  e^{-t/x - (1-\lambda)x}  x^{\nu+j} dx =   \sum_{j=0}^n S(n,j)\  {\lambda^j  \rho_{\nu+j+1}\left( t(1-\lambda) \right) \over (1-\lambda)^{\nu+j+1}}  $$
which yields (2.3).  To calculate moments (2.2) we derive

$$\gamma_n =   \int_0^\infty  e^{-t/x -x}  x^{\nu} \sum_{k=0}^\infty { (x\lambda)^k  \over (k-n)!} =  \lambda^n \int_0^\infty  e^{-t/x - (1-\lambda) x}  x^{\nu+n} dx$$

$$=  { \lambda^n   \rho_{\nu+n+1}\left( t(1-\lambda) \right) \over (1-\lambda)^{\nu+n+1}}.$$
This completes the proof of Lemma 1. 

\end{proof}

Now we will interpret  discrete orthogonality conditions (1.12), (1.13) in terms of the so-called composition orthogonality.

{\bf Theorem 1}. {\it Let $n \in \mathbb{N},\ t > 0,\ 0 < \lambda < 1,\ \nu > -1$. Orthogonality conditions $(1.12), (1.13)$  are equivalent to the composition orthogonality relatively to the function $e^{\lambda x}$ with respect to the measure $e^{-t/x -x}  x^{\nu} dx$.   Precisely, the following equalities hold, respectively, }

$$  \int_0^\infty  e^{-t/x -x}  x^{\nu}   P_n\left( x D; t,\lambda \right) \left( x D \right)^j \left\{ e^{\lambda x} \right\} dx = 0,\quad   j =0,\dots, n-1,\eqno(2.7)$$

$$  \int_0^\infty  e^{-t/x -x}  x^{\nu}    P_n\left( x D; t,\lambda\right) \left\{  x^j e^{\lambda x} \right\} dx = 0,\quad   j =0,\dots, n-1.\eqno(2.8)$$

\begin{proof} Recalling (1.2), we rewrite (1.12)  in the form

$$ \int_0^\infty  e^{-t/x -x}  x^{\nu} \sum_{k=0}^\infty P_n(k; t,\lambda)  { (x\lambda)^k k^j\over k!}  dx= 0,\ j =0,\dots, n-1,\eqno(2.9)$$ 
where the interchange of the order of integration and summation is justified as above by virtue of the dominated convergence theorem.   Then proceeding as in (2.5), we obtain

$$ \int_0^\infty  e^{-t/x -x}  x^{\nu} \sum_{k=0}^\infty P_n(k; t,\lambda)  { (x\lambda)^k k^j\over k!}  dx =  \int_0^\infty  e^{-t/x -x}  x^{\nu}  \left( x D \right)^j \sum_{k=0}^\infty P_n(k; t,\lambda)  { (x\lambda)^k \over k!}  dx$$

$$=  \int_0^\infty  e^{-t/x -x}  x^{\nu}  \left( x D \right)^j  P_n\left( x D; t,\lambda\right) \left\{ e^{\lambda x} \right\} dx =  \int_0^\infty  e^{-t/x -x}  x^{\nu}   P_n\left( x D; t,\lambda\right) \left( x D \right)^j \left\{ e^{\lambda x} \right\} dx.$$
Consequently, we derive  equalities  (2.7).  In order to prove (2.8) we have, similarly,  from (1.13)

$$\sum_{k=0}^\infty P_n(k; t,\lambda)  \rho_{k+\nu+1}(t) { \lambda^k (-k)_j  \over k!} =  (-1)^j \int_0^\infty  e^{-t/x -x}  x^{\nu}  \sum_{k=0}^\infty P_n(k; t,\lambda)  { (x\lambda)^k  \over (k-j)!} dx$$ 

$$=  (-1)^j \int_0^\infty  e^{-t/x -x}  x^{\nu}    P_n\left( x D; t,\lambda\right)  \left\{ \sum_{k=0}^\infty { (x\lambda)^k  \over (k-j)!} \right\} dx =   (-1)^j \lambda^j \int_0^\infty  e^{-t/x -x}  x^{\nu}    P_n\left( x D; t,\lambda\right) \left\{  x^j e^{\lambda x} \right\} dx.$$
Hence we arrive at the conditions (2.8).

\end{proof}

{\bf Remark 1}. {\it Orthogonal polynomials for the weight $e^{-t/x -x}  x^{\nu} $ were  investigated recently by the author in \cite{YaK}. }

{\bf Corollary 1}. {\it  Orthogonality relations  $(1.12), (1.13)$ are equivalent to the conditions

$$ \int_0^\infty   e^{-(1-\lambda) x- t/x}   f_n(x; t,\lambda)\  x^{\nu+j}  dx =0,\quad   j =0,\dots, n-1,\eqno(2.10)$$
where 

$$f_n(x; t,\lambda)= \sum_{m=0}^n  \sum_{k=0}^m  \sum_{i=0}^k  (-1)^{k+m}  a_{n,m}  S(m+1,k+1) \  {x^i \ k! \over i!}  L_{k-i}^{-\nu-1} \left({t\over x}\right)\eqno(2.11)$$
and $a_{n,m}$ are coefficients of $P_n(x; t,\lambda)$.}

\begin{proof} Indeed, taking (2.7), we integrate by parts, eliminating the integrated terms.  Then we find

$$  \int_0^\infty  e^{-t/x -x}  x^{\nu}   P_n\left( x D; t,\lambda\right) \left( x D \right)^j \left\{ e^{\lambda x} \right\} dx $$

$$=  \int_0^\infty  P_n\left( - D x; t,\lambda \right)  \left\{  e^{-t/x -x}  x^{\nu} \right\}  \left( x D \right)^j \left\{ e^{\lambda x} \right\} dx.\eqno(2.12) $$
Hence, employing a companion of the equality (2.6) from \cite{Riordan}

$$ \left(  D  x\right)^m = \sum_{k=0}^m S(m+1,k+1) x^k D^k,\eqno(2.13)$$
we calculate the composition $\left(D x \right)^m   \left\{  e^{-t/x -x}  x^{\nu} \right\}$, invoking  Entry 1.1.3.2 on p. 4 in \cite{Bry}.  Therefore we obtain

$$ \left(D x \right)^m   \left\{  e^{-t/x -x}  x^{\nu} \right\} = \sum_{k=0}^m S(m+1,k+1) \ x^k D^k   \left\{  e^{-t/x -x}  x^{\nu} \right\}$$

$$=  e^{-x} \sum_{k=0}^m S(m+1,k+1) \ x^k   \sum_{i=0}^k (-1)^i \binom{k}{i}  D^{k-i} \left\{  e^{-t/x}  x^{\nu} \right\} $$

$$=  e^{-x- t/x} x^\nu \sum_{k=0}^m  \sum_{i=0}^k  (-1)^k S(m+1,k+1) \  {x^i \ k! \over i!}  L_{k-i}^{-\nu-1} \left({t\over x}\right).$$
Moreover, via (2.6) we write

$$ \left( x D \right)^j \left\{ e^{ax} \right\}  =  \sum_{i=0}^j S(j,i) x^i D^i \left\{ e^{\lambda x} \right\} =  e^{\lambda x}  \sum_{i=0}^j S(j,i) (\lambda x)^i.$$
Thus, substituting these expressions  in (2.12), where 

$$P_n(x; t,\lambda)= \sum_{m=0}^n a_{n,m}  x^m,\ a_{n,n} \equiv a_n,\ a_{n,n-1} \equiv b_n,$$
it gives 

$$ \int_0^\infty   e^{-(1-\lambda ) x- t/x}   f_n(x; t,\lambda)\  x^{\nu} \left(\sum_{i=0}^j S(j,i) (\lambda x)^i\right) dx =0,\quad   j =0,\dots, n-1,\eqno(2.14)$$
where 

$$f_n(x; t,\lambda)= \sum_{m=0}^n  \sum_{k=0}^m  \sum_{i=0}^k  (-1)^{k+m}  a_{n,m}  S(m+1,k+1) \  {x^i \ k! \over i!}  L_{k-i}^{-\nu-1} \left({t\over x}\right).$$
But (2.14) is equivalent to (2.10). In the same manner conditions (1.13) can be treated. 

\end{proof}

The difference operators $\Delta$ and $\nabla$ are defined accordingly

$$\Delta f = f(x+1)- f(x), \quad \nabla f = f(x)- f(x-1).\eqno(2.15)$$
It is used to establish the structural relation for the orthogonal polynomials $P_n(x; t,\lambda)$. But first we observe via (1.7), (1.8) that these polynomials are Charlier polynomials \cite{Sze} in the sense of composition orthogonality, i.e.  (1.1) reads

$$\sum_{k=0}^\infty P_n(k; t,\lambda) P_m(k; t,\lambda)  { (\lambda I_-)^k\over k!}  \left\{ \rho_{\nu+1}(t) \right\} =  \delta_{n,m},\quad 0 < \lambda  < 1.\eqno(2.16)$$ 
Now, expanding $P_n(x+1; t,\lambda)$ in terms of $P_n$, we have

$$P_n(x+1; t,\lambda) = \sum_{j=0}^n c_{n,j}  P_j(x; t,\lambda),\eqno(2.17)$$
where, evidently, $c_{n,n} =1$.   Other coefficients are given by the formula (see (1.1))

$$ c_{n,j} =  \sum_{k=0}^\infty P_n(k+1; t,\lambda) P_j(k; t,\lambda) \omega_k.\eqno(2.18)$$
Then by virtue of (1.18) we deduce from (2.18)

$$ c_{n,j} =  - {1\over \lambda } \sum_{k=0}^\infty P_n(k; t,\lambda)\ k P_j(k-1; t,\lambda) {\partial \omega_k\over \partial t} =  - {1\over \lambda} {\partial \over \partial t} \left\{ \sum_{k=0}^\infty P_n(k; t,\lambda)\  k P_j(k-1; t,\lambda) \omega_k  \right\}$$

$$+  {1\over \lambda} \sum_{k=0}^\infty {\partial P_n\over \partial t} (k; t,\lambda)\  k P_j(k-1; t,\lambda) \omega_k (t,\lambda) +  {1\over \lambda} \sum_{k=0}^\infty P_n(k; t,\lambda)\  k {\partial P_j\over \partial t} (k-1; t,\lambda) \omega_k ,\eqno(2.19)$$ 
where we assume that polynomial coefficients are $C^1(\mathbb{R}^2_+)$-functions of variables $t, \lambda$  and the differentiation under the series  sign is allowed on closed intervals,  owing to the estimate

$$  \sum_{k=0}^\infty \left|  {\partial \over \partial t} \left\{ P_n(k; t,\lambda)\  k P_j(k-1; t,\lambda) \omega_k  \right\} \right|$$

$$ \le \sum_{m=0}^n \sum_{i=0}^j  \max \left\{ \max_{t \in [\alpha, \beta]}   \left|  {\partial \over \partial t} \left\{ a_{n,m} (t,\lambda)  a_{j,i} (t,\lambda) \right\} \right|,  
 \max_{t \in [\alpha, \beta]}   \left|  a_{n,m} (t,\lambda)  a_{j,i} (t,\lambda)  \right| \right\} $$

$$\times   \sum_{k=1}^\infty  k^{m+1} (k-1)^i  (k+\nu+1) \Gamma (k+\nu) {\lambda^k\over k!} < \infty, \quad 0 < \lambda < 1,\ \nu >  -1.$$
Hence by orthogonality we get

$$ c_{n,j} =   {1\over \lambda} \sum_{k=0}^\infty {\partial P_n\over \partial t} (k; t,\lambda)\  k P_j(k-1; t,\lambda) \omega_k ,\quad j <  n-1.\eqno(2.20)$$
When $j= n-1$, we have from (1.1), (1.14)

$$ \sum_{k=0}^\infty P_n(k; t,\lambda)\  k P_{n-1}(k-1; t,\lambda) \omega_k = A_n,\eqno(2.21)$$

$$\sum_{k=0}^\infty P_n(k; t,\lambda)\  k {\partial P_{n-1}\over \partial t} (k-1; t,\lambda) \omega_k  = {1\over a_n} {\partial a_{n-1}\over  \partial t},\eqno(2.22)$$
and therefore we end up from (2.19) with the value 

$$c_{n,n-1}=  {1\over \lambda }  { A_n \over a_n}  {\partial a_{n}\over \partial t}  +  {1\over \lambda } \sum_{k=0}^\infty {\partial P_n\over \partial t} (k; t,\lambda)\  k P_{n-1}(k-1; t,\lambda) \omega_k .\eqno(2.23)$$
In the meantime, recalling (2.20) and the Christoffel-Darboux formula (1.16), we calculate the sum

$$\sum_{j=0}^{n-2} c_{n,j} P_j(x; t,\lambda) =   {1\over \lambda} \sum_{k=0}^\infty {\partial P_n\over \partial t} (k; t,\lambda)\  k \omega_k   \sum_{j=0}^{n-2}  P_j(x; t,\lambda) P_j(k-1; t,\lambda) $$

$$=  {A_{n-1}  \over \lambda}   \sum_{k=0}^\infty   {\partial P_n\over \partial t} (k; t,\lambda)  \frac{  k \omega_k  }{x-k+1} \bigg[  P_{n-1}(x; t,\lambda) P_{n-2}(k-1; t,\lambda) -  P_{n-2}(x; t,\lambda) P_{n-1}(k-1; t,\lambda)\bigg].\eqno(2.24) $$
Thus, combining with (2.17), (2.19), (2.21), (2.22), (2.24) we establish  the structural relation for the sequence $\{P_n(x; t,\lambda)\}_{n\ge 0}$.

{\bf Theorem 2}. {\it For the orthogonal polynomials $P_n(x; t,\lambda)$ it has the equality }

$$  P_n (x+1;  t,\lambda) =  P_n(x; t,\lambda)+  { A_n \over \lambda a_n }\  {\partial a_n\over \partial t}    P_{n-1}(x; t,\lambda) $$

$$+  {1\over \lambda}  \sum_{k=0}^\infty  {\partial P_n\over \partial t} (k; t,\lambda)  \frac{k \omega_k  }{x-k+1} \bigg[  P_{n-1}(x; t,\lambda) \bigg[ A_{n-1}  P_{n-2}(k-1; t,\lambda) + (x-k+1)  P_{n-1}(k-1; t,\lambda)\bigg] $$

$$\bigg. -  A_{n-1} P_{n-2}(x; t,\lambda) P_{n-1}(k-1; t,\lambda)\bigg].\eqno(2.25)$$

{\bf Corollary 2}. {\it  The following identity holds

$$ \sum_{k=0}^\infty {\partial P_n\over \partial t} (k; t,\lambda)  \  k P_{n-1} (k-1; t,\lambda) \omega_k   =  { \lambda n\over A_n}  -  {A_n\over a_n}\  {\partial a_{n}\over  \partial t} .\eqno(2.26)$$
Besides, $c_{n,n-1}= n/ A_n $}.

\begin{proof}  Writing (2.25) in terms of the forward difference $\Delta$ (2.15), we have 

$$  \Delta P_n (x;  t,\lambda) =  { A_n \over \lambda a_n}  {\partial a_{n}\over  \partial t} P_{n-1}(x; t,\lambda) + {1\over \lambda}  P_{n-1}(x; t,\lambda) \sum_{k=0}^\infty  {\partial P_n\over \partial t} (k; t,\lambda) \  k  P_{n-1}(k-1; t,\lambda) \omega_k $$

$$+  {1\over \lambda}   A_{n-1} \sum_{k=0}^\infty   {\partial P_n\over \partial t} (k; t,\lambda)  \frac{k \omega_k  }{x-k+1} $$

$$\times \bigg[  P_{n-1}(x; t,\lambda) P_{n-2}(k-1; t,\lambda)  - P_{n-2}(x; t,\lambda) P_{n-1}(k-1; t,\lambda)\bigg].\eqno(2.27)$$
Then, comparing the leading coefficient in the latter equality, we get (2.26). The value of $c_{n,n-1}$ follows immediately from (2.23).

\end{proof}

{\bf Corollary 3.} {\it In terms of the forward difference  equality $(2.25)$ reads }

$$\Delta P_n (x; t,\lambda) =   {n\over  A_n}\   P_{n-1}(x; t,\lambda) +  {A_{n-1} \over \lambda}  \sum_{k=0}^\infty  {\partial P_n\over \partial t} (k; t,\lambda) \frac{ k \omega_k  }{x-k+1}$$

$$\times  \bigg[  P_{n-1}(x; t,\lambda) P_{n-2}(k-1; t,\lambda)  -  P_{n-2}(x; t,\lambda) P_{n-1}(k-1; t,\lambda)\bigg].\eqno(2.28)$$
Further, if we differentiate (1.1) with $n=m$, it gives

$$2  \sum_{k=0}^\infty P_n(k; t,\lambda)  {\partial P_n\over \partial t} (k; t,\lambda) \omega_k +  \sum_{k=0}^\infty P^2_n(k; t,\lambda) {\partial \omega_k\over \partial t} = 0.$$
Hence

$$ \sum_{k=0}^\infty P^2_n(k; t,\lambda) {\partial \omega_k\over \partial t}  = -   {2 \over a_n(t)} {\partial a_n\over \partial t} .\eqno(2.29)$$
However, the use of the differential equation (1.18) allows us to compute the sum $\sum_{k=0}^\infty P^2_n(k+1,t)\omega_k.$  In fact, we have

$$ \sum_{k=0}^\infty P^2_n(k+1; t,\lambda)\omega_k = \sum_{k=1}^\infty P^2_n(k; t,\lambda)\omega_{k-1} = - {1\over \lambda} \sum_{k=0}^\infty k P^2_n(k; t,\lambda) {\partial \omega_k \over \partial t}$$

$$= - {1\over \lambda}  {\partial  \over \partial t} \left( \sum_{k=0}^\infty k P^2_n(k; t,\lambda)\omega_k\right) +  {2\over \lambda} \sum_{k=0}^\infty k P_n(k; t,\lambda) {\partial P_n\over \partial t} (k; t,\lambda)\omega_k.$$
Hence (1.1), (1.14), (1.15) imply 

$$ \sum_{k=0}^\infty P^2_n(k+1; t,\lambda)\omega_k = - {1\over \lambda} {\partial B_n\over \partial t}  + 
 {2\over \lambda}\bigg[ {1\over a_n}  {\partial b_n\over \partial t}  - { b_{n+1}\over a_{n} a_{n+1}}  {\partial a_n\over \partial t}\bigg],$$
i.e. after slight simplification we find the equality

$$ \sum_{k=0}^\infty P^2_n(k+1,t)\omega_k = {1\over \lambda} \bigg[  {\partial \over \partial t} \left({b_{n} \over a_{n}} + {b_{n+1}\over a_{n+1}} \right)   -   {2 B_n\over a_n}  {\partial a_n\over \partial t}\bigg].\eqno(2.30)$$
On the other hand, we have from (2.24), (2.28),  (1.1)

$$ \sum_{k=0}^\infty P^2_n(k+1,t)\omega_k = 1+  \ {n^2\over A^2_n} + \sum_{j=0}^{n-2} c^2_{n,j}.\eqno(2.31)$$
Thus, combining with (2.30), we derive

$$  \sum_{j=0}^{n-2} c^2_{n,j} =  {1\over \lambda} \bigg[ {\partial \over \partial t} \left({b_{n} \over a_{n}} + {b_{n+1}\over a_{n+1}} \right) -    {2B_n\over a_n} {\partial a_n\over \partial t}\bigg]- 1 -   {n^2\over A^2_n},\eqno(2.32)$$
where $c_{n,j},\ j < n-1$ is defined by (2.20).

Further,  decomposing  $P_n(x+2; t,\lambda)$ in terms of $P_n$, we find

$$P_n(x+2; t,\lambda) = \sum_{j=0}^n d_{n,j}(t,\lambda) P_j(x; t,\lambda).\eqno(2.33)$$
Hence $d_{n,n} =1$ and 

$$ d_{n,j} =  \sum_{k=0}^\infty P_n(k+2; t,\lambda) P_j(k; t,\lambda) \omega_k.\eqno(2.34)$$
Then  (1.17), (1.18)  suggest  the equalities 

$$ d_{n,j} =   \sum_{k=2}^\infty P_n(k; t,\lambda)\  P_j(k-2; t,\lambda) \omega_{k-2} =  {1\over t\lambda^2}  \sum_{k=0}^\infty P_n(k; t,\lambda)\  k(k-1) P_j(k-2; t,\lambda) \omega_{k}  $$

$$- {1\over t\lambda}  \sum_{k=1}^\infty P_n(k; t,\lambda)\  (k-1) (k+\nu) P_j(k-2; t,\lambda) \omega_{k-1}  =  {1\over t\lambda^2}  \sum_{k=0}^\infty P_n(k; t,\lambda)\  k(k-1) P_j(k-2; t,\lambda) \omega_{k} $$

$$ + {1\over t\lambda^2} {\partial \over \partial t}  \left(\sum_{k=0}^\infty P_n(k; t,\lambda)\  k(k-1) (k+\nu) P_j(k-2; t,\lambda) \omega_{k} \right)$$

$$-  {1\over t\lambda^2}  \sum_{k=0}^\infty {\partial P_n\over \partial t} (k; t,\lambda)\  k(k-1) (k+\nu) P_j(k-2; t,\lambda) \omega_{k}  -  {1\over t\lambda^2}  \sum_{k=0}^\infty P_n(k,t)\  k(k-1) (k+\nu) $$

$$\times {\partial P_j\over \partial t} (k-2; t,\lambda) \omega_{k} .\eqno(2.35) $$
Therefore the same analysis shows

$$d_{n,j} =  -  {1\over t\lambda^2}  \sum_{k=0}^\infty {\partial P_n\over \partial t} (k; t,\lambda)\  k(k-1) (k+\nu) P_j(k-2; t,\lambda) \omega_{k} ,\quad j < n-3.\eqno(2.36)$$
Then,  accordingly, we have from (2.35), (1.1), (1.14), (1.15) 

$$ d_{n,n-3} = {1\over t\lambda^2}  {\partial \over \partial t} \left(\sum_{k=0}^\infty P_n(k; t,\lambda)\  k(k-1) (k+\nu) P_{n-3}(k-2; t,\lambda) \omega_{k} \right)$$

$$-  {1\over t\lambda^2}  \sum_{k=0}^\infty {\partial P_n\over \partial t} (k; t,\lambda)\  k(k-1) (k+\nu) P_{n-3}(k-2; t,\lambda) \omega_{k} $$

$$ -  {1\over t\lambda^2}  \sum_{k=0}^\infty P_n(k; t,\lambda)\  k(k-1) (k+\nu)  {\partial P_{n-3}\over \partial t} (k-2; t,\lambda) \omega_{k} $$

$$= -  {1\over t\lambda^2} \bigg[  {a_{n-3}  \over a^2_n } {\partial a_n\over \partial t}  +  \sum_{k=0}^\infty  {\partial P_n\over \partial t} (k; t,\lambda)\  k(k-1) (k+\nu)  P_{n-3}(k-2; t,\lambda) \omega_{k} \bigg],\eqno(2.37)$$

$$ d_{n,n-2} =   {1\over t\lambda^2} \bigg[ {a_{n-2}\over a_n} - \bigg[ (\nu+2n-5) a_{n-2} + b_{n-2}\bigg] \ {1 \over a^2_n} {\partial a_n\over \partial t} -   a_{n-2}{\partial \over \partial t} \left( {b_{n+1} \over a_{n+1} a_n} \right) \bigg.$$

$$\bigg.  -    \sum_{k=0}^\infty  {\partial P_n\over \partial t} (k; t,\lambda)\  k(k-1) (k+\nu) P_{n-2}(k-2; t,\lambda) \omega_{k}  \bigg].\eqno(2.38) $$
Now, writing $P_n$ as 

$$P_n(x; t,\lambda)= a_n x^n+ b_n x^{n-1} + c_n x^{n-2} + \hbox{lower degrees},\eqno(2.39)$$
we recall (1.1), (1.12) to get  in a straightforward way the equality

$$ \sum_{k=0}^\infty P_n(k; t,\lambda)  k^{n+2}  \omega_{k}  = {b_{n+2} b_{n+1} \over a_{n+2}  a_{n+1} a_n} -  {c_{n+2} \over a_{n+2}  a_n}.\eqno(2.40)$$ 
This identity allows to obtain the value of the coefficient $d_{n,n-1}$.  In fact, we derive 

$$ d_{n,n-1} =   {1\over t\lambda^2}  \bigg[ A_n \left[ B_{n}   + B_{n-1} +1-2n\right]  \bigg.$$

$$+ a_{n-1}  {\partial \over \partial t} \bigg(   {b_{n+2} b_{n+1}  \over a_{n+2}  a_{n+1} a_n} -  {c_{n+2}  \over a_{n+2}  a_n}\bigg)$$

$$- \bigg[ b_{n-1} + (\nu+1-2n) a_{n-1} \bigg]  {\partial \over \partial t} \bigg( { b_{n+1}  \over  a_{n+1} a_n} \bigg)  $$

$$ + \bigg[  (\nu+2(n-1))  a_{n-1} -   (\nu-1+ 2(n-2)^2 ) b_{n-1} -  c_{n-1}\bigg]  { 1  \over  a^2_n} {\partial a_n\over \partial t}$$

$$\bigg. -   \sum_{k=0}^\infty {\partial P_n\over \partial t} (k; t,\lambda)\  k(k-1) (k+\nu) P_{n-1}(k-2; t,\lambda) \omega_{k} \bigg].\eqno(2.41)$$
Then, writing (2.33) in the form

$$P_n(x+2; t,\lambda) - P_n(x; t,\lambda) =  d_{n,n-1}  P_{n-1}(x; t,\lambda) +  d_{n,n-2}  P_{n-2}(x; t,\lambda)$$

$$+  d_{n,n-3}  P_{n-3}(x; t,\lambda)+ \sum_{j=0}^{n-4} d_{n,j} P_j(x; t,\lambda),\eqno(2.42)$$
we compare leading coefficients to get  $d_{n,n-1}= 2n/A_n$. Hence from (2.41) we find 

$$ {1\over t\lambda^2} \sum_{k=0}^\infty  {\partial P_n\over \partial t}  (k; t,\lambda)\  k(k-1) (k+\nu) P_{n-1}(k-2; t,\lambda) \omega_{k}  $$

$$=   {1\over t\lambda^2}  \bigg[ A_n \left[ B_{n}  + B_{n-1} +1-2n\right] +  a_{n-1} {\partial \over \partial t}  \bigg(   {b_{n+2} b_{n+1}  \over a_{n+2}  a_{n+1} a_n} -  {c_{n+2}  \over a_{n+2}  a_n}\bigg)\bigg.$$

$$- \bigg[ b_{n-1} + (\nu+1-2n) a_{n-1} \bigg]  {\partial \over \partial t} \bigg( { b_{n+1}  \over  a_{n+1} a_n} \bigg)  $$

$$ + \bigg[  (\nu+ 2(n-1)) a_{n-1} -   (\nu-1+2(n-2)^2) b_{n-1} -  c_{n-1}\bigg]  { 1 \over  a^2_n} {\partial a_n\over \partial t}  \bigg] - {2n\over A_n} .$$
Furthermore, in the same manner as above by virtue of the Christoffel-Darboux formula (1.16) and identities (2.37), (2.38), (2.41) we establish the following theorem.

{\bf Theorem 3}. {\it For the orthogonal polynomials $P_n(x; t,\lambda)$ it has the equality }

$$P_n(x+2; t,\lambda) =  P_n(x; t,\lambda) + {2n\over A_n(t,\lambda)} P_{n-1} (x; t,\lambda)  -  {1\over t\lambda^2}  {a_{n-3}  \over a^2_n} {\partial a_n\over \partial t}  P_{n-3} (x; t,\lambda)  $$

$$+  {1\over t\lambda^2} \bigg[ {a_{n-2}\over a_n} - \bigg[ (\nu+2n-5) a_{n-2}+ b_{n-2} \bigg] \ {1 \over a^2_n} {\partial a_n\over \partial t}   -   a_{n-2} {\partial \over \partial t}  \left( {b_{n+1}\over a_{n+1} a_n} \right) \bigg] P_{n-2} (x; t,\lambda) $$

$$- {A_{n-1} \over t\lambda^2} \sum_{k=0}^\infty  {\partial P_n\over \partial t}  (k; t,\lambda) \frac{ k(k-1) (k+\nu) \omega_k } {x-k+2}$$

$$\times  \bigg[  P_{n-1}(k; t,\lambda) P_{n-2}(k-2; t,\lambda)  -  P_{n-2}(x; t,\lambda) P_{n-1}(k-2; t,\lambda)\bigg].\eqno(2.43)$$

\section{Recurrence  relations}

Let us differentiate (1.12) with respect to $\lambda$ under the same justification as in (2.19). Then, employing (1.1), (1.14) and (1.19), we obtain

$$\sum_{k=0}^\infty  \bigg[ \lambda {\partial P_n\over \partial \lambda} (k; t,\lambda)  + A_n(t,\lambda) P_{n-1} (k; t,\lambda) \bigg] \omega_k k^j = 0,\ j =0,\dots, n-1.\eqno(3.1)$$ 
The latter orthogonality conditions immediately imply the differential-difference equation

$$ \lambda {\partial P_n\over \partial \lambda} (x; t,\lambda)  + A_n P_{n-1} (x; t,\lambda) = C_n P_n(x; t,\lambda).\eqno(3.2)$$
The constant $ C_n$ can be obtained, comparing leading coefficients in (3.2). Hence it yields 

$$ \lambda  a_n {\partial P_n\over \partial \lambda} (x; t,\lambda)  + a_{n-1}  P_{n-1} (x; t,\lambda) = \lambda {\partial a_n\over \partial \lambda}\  P_n(x; t,\lambda).\eqno(3.3)$$
On the other hand, differentiating by $\lambda$ equality (1.1) when $m=n$, we find

$$ 2 \sum_{k=0}^\infty P_n(k; t,\lambda) {\partial P_n\over \partial \lambda} (k; t,\lambda) \omega_k + \sum_{k=0}^\infty P^2_n(k; t,\lambda) {\partial  \omega_k \over \partial \lambda}  = 0.$$
Hence, recalling (1.14), (1.19), we get the identity

$$B_n = - 2 {\lambda \over a_n} {\partial a_n\over \partial \lambda} .\eqno(3.4)$$
Moreover, relation (3.4) can be rewritten in the form

$$ \lambda  {\partial P_n\over \partial \lambda} (x; t,\lambda)  + A_n P_{n-1} (x; t,\lambda) = - {B_n \over 2} P_n(x; t,\lambda).\eqno(3.5)$$
But the differentiation of $B_n$ by $\lambda$ yields, in turn,

$$ {\partial B_n\over \partial \lambda} = 2 \sum_{k=0}^\infty k P_n(k; t,\lambda) {\partial P_n  \over \partial \lambda} (k; t,\lambda) \omega_k +  {1\over \lambda} \sum_{k=0}^\infty \bigg[ k P_n(k; t,\lambda)\bigg]^2  \omega_k .$$ 
Hence via orthogonality and (3.5) we end up after simplification with the identity

$$A^2_{n+1} + A^2_n+   \lambda    {\partial \over \partial \lambda} \left(  {b_{n} \over a_{n} }  + {b_{n+1} \over a_{n+1}} \right) = 0.\eqno(3.6)$$
However, equating  coefficients in front of $x^{n-1} $ in (3.5), we deduce a more simple relation 

$$A^2_n +   \lambda    {\partial \over \partial \lambda} \left(  {b_{n} \over a_{n} } \right) = 0.\eqno(3.7)$$
It can be rewritten in the form 

$$\lambda    {\partial B_n \over \partial \lambda} =  A^2_{n+1}  - A^2_n.\eqno(3.8)$$
It is not difficult to verify the equation

$$\lambda    {\partial A^2_n \over \partial \lambda} =   A^2_n \bigg(B_n - B_{n-1}\bigg).\eqno(3.9)$$
Equations (3.8), (3.9) constitute the classical Toda system by $\lambda$ for recurrence coefficients (1.15).  Meanwhile,  the orthogonality relation implies

$$ \sum_{k=0}^\infty   {\partial P_n  \over \partial \lambda} (k; t,\lambda) P_{n-1} (k; t,\lambda) \omega_k  + \sum_{k=0}^\infty   P_n (k; t,\lambda) P_{n-1} (k; t,\lambda) {\partial    \omega_k \over \partial \lambda} = 0,$$
i.e. via (1.19), (3.4) we get the identity

$$ a_n b_n B_n + 2 a_{n-1}^2  = 0.\eqno(3.10)$$
Moreover,  as a consequence of  (3.7), (3.8), (3.9) the latter equality yields

$$ B_{n-1}  B_n B_{n+1}  =  A^2_{n+1}  \bigg(B_n - B_{n+1}\bigg).\eqno(3.11)$$
In fact, we have from (3.10)

$$ {b_{n} \over a_{n}}  B_n + 2 A^2_n= 0.\eqno(3.12)$$
After differentiation and multiplication by $\lambda$ of the latter equality and the use of (3.7), (3.8), (3.9) we find

$$ A^2_n B_n - A^2_n B_{n-1}  + {b_{n} \over a_{n} }  A^2_{n+1}  - A^2_n  {b_{n-1} \over a_{n-1} }   =0.$$
Hence

$${b_{n-1} \over a_{n-1} }- {b_{n}\over a_{n} } = {1\over 2} \left[  {b_{n} \over a_{n}}  {A^2_{n+1}\over A^2_{n} } - {b_{n+1}\over a_{n+1} } \right].$$
Meanwhile, the use of (3.12) implies

$$B_{n-1}  =  A^2_{n+1} \bigg[ {1\over B_{n+1}} -   {1\over B_{n}}\bigg],$$
which, in turn, yields (3.11).

Now, recalling  the orthogonality conditions (1.13), we make summation by parts to deduce the equalities

$$\sum_{k=0}^\infty \Delta \bigg[ P_n(k; t,\lambda)  \rho_{k+\nu+1}(t) \bigg] { \lambda^k \over (k-j )!} = 0,\quad j =0,\dots, n-2.\eqno(3.13)$$ 
We will need to take into account the dependence of the polynomials $P_n$ upon $\nu$, i.e. we write $P_n\equiv P^\nu_n$. Then (3.13) is equivalent to the conditions

$$\sum_{k=0}^\infty  P^\nu_n(k+1; t,\lambda)  \rho_{k+\nu+2}(t)  { \lambda^k \over (k-j )!} = 0,\quad j =0,\dots, n-2.\eqno(3.14)$$ 
This means that the sequence $\{P_n^\nu(x+1; t, \lambda) \}_{n\ge 0}$ is quasi-orthogonal with respect to the weight $w_k^{\nu+1}$. Hence, expanding $P_n^\nu(x+1)$ in terms of the sequence   $\{P_n^{\nu+1}\}_{n\ge 0}$

$$P_n^\nu\left(x+1; \ t,\lambda\right)= \sum_{m=0}^n \gamma^\nu_{n,m}  P_m^{\nu+1}\left(x; \ t,\lambda\right),$$
one proves  in a straightforward way, owing to (3.14),  that $\gamma^\nu_{n,m} = 0,\ m=0,\dots, n-2$. Then, accordingly,

$$ P_n^\nu\left(x+1; \ t,\lambda\right)= \gamma^\nu_{n,n}  P_n^{\nu+1}\left(x; \ t,\lambda\right) + \gamma^\nu_{n,n-1}  P_{n-1}^{\nu+1}\left(x; \ t,\lambda\right),\eqno(3.15)$$
where

$$ \gamma^\nu_{n,n-1}  = \sum_{k=0}^\infty  P^\nu_n(k+1; t,\lambda)  P_{n-1}^{\nu+1}\left(k; \ t,\lambda\right) \rho_{k+\nu+2}(t)  { \lambda^k \over k!},\eqno(3.16)$$

$$ \gamma^\nu_{n,n}  = \sum_{k=0}^\infty  P^\nu_n(k+1; t,\lambda)  P_{n}^{\nu+1}\left(k; \ t,\lambda\right) \rho_{k+\nu+2}(t)  { \lambda^k \over k!}.\eqno(3.17)$$
Then, recalling  (2.15), we  write the structural relation (2.28) for polynomials $P_n^\nu$ in the form

$$ \Delta P_n^\nu\left(x; t,\lambda\right)= \gamma^\nu_{n,n} \nabla P_n^{\nu+1}\left(x; \ t,\lambda\right) + \gamma^\nu_{n,n-1}  \nabla P_{n-1}^{\nu+1}\left(x; \ t,\lambda\right).\eqno(3.18)$$
Hence from the orthogonality (1.1) it is easily seen the values

$$ \gamma^\nu_{n,n-1}  =   { a_n^\nu \over a_{n-1}^{\nu+1} } \bigg[ n + {b_n^\nu \over a_{n}^{\nu}} - {b_n^{\nu+1}  \over a_n^{\nu+1} } \bigg],\eqno(3.19)$$

$$ \gamma^\nu_{n,n} = {a_n^\nu \over a_n^{\nu+1} },\eqno(3.20)$$
$$ \gamma^\nu_{n,n-1}  = { \gamma^\nu_{n,n}  \over A_n^{\nu+1} } \bigg[ n + {b_n^\nu \over a_{n}^{\nu} } - {b_n^{\nu+1}  \over a_n^{\nu+1}  } \bigg],\eqno(3.21)$$
On the other hand,  equality (3.16) reads via (1.1), (1.14), (1.15)

$$ \gamma^\nu_{n,n-1}  = {1\over \lambda} \sum_{k=1}^\infty  P^\nu_n(k; t,\lambda)  P_{n-1}^{\nu+1}\left(k-1; \ t,\lambda\right) \rho_{k+\nu+1}(t)  { \lambda^k \over (k-1)!}$$

$$= {A_n^\nu \over \lambda} \sum_{k=0}^\infty  P^\nu_{n-1} (k; t,\lambda)  P_{n-1}^{\nu+1}\left(k-1; \ t,\lambda\right) \rho_{k+\nu+1}(t)  { \lambda^k \over k!}$$

$$= {1\over \lambda} {a_{n-1}^{\nu+1} \over a_{n}^\nu }=   {1\over \lambda} {A_n^{\nu} \over   \gamma^\nu_{n-1,n-1} } = {1\over \lambda} {A_n^{\nu+1} \over   \gamma^\nu_{n,n} }.\eqno(3.22) $$
Thus, combining with (3.19), we derive the following recurrence relations

$$\bigg[  {a_{n-1}^{\nu+1}\over a_{n}^\nu }\bigg]^2 = \bigg[{A_n^{\nu+1} \over   \gamma^\nu_{n,n} }\bigg]^2=  \lambda \bigg[ n + {b_n^\nu \over a_{n}^{\nu} } - {b_n^{\nu+1} \over a_n^{\nu+1}  } \bigg].\eqno(3.23)$$
Analogously, we treat equality (3.17).  It gives,

$$  \gamma^\nu_{n,n}  = {1\over \lambda } \sum_{k=0}^\infty  P^\nu_n(k; t,\lambda)  k P_{n}^{\nu+1}\left(k-1; \ t,\lambda\right) \rho_{k+\nu+1}(t)  { \lambda^k \over k!}$$

$$= {1\over \lambda }  {a_n^{\nu+1} \over a_n^{\nu}}  \bigg[ {b_n^{\nu+1} \over a_{n}^{\nu+1} } -  { b_{n+1}^\nu\over  a_{n+1}^{\nu}} - n \bigg].$$
Hence, recalling (3.20),  we deduce the identity

$$  \left[\gamma^\nu_{n,n} \right]^2  = \bigg[ {a_n^\nu \over a_n^{\nu+1} } \bigg]^2 =  {1\over \lambda } \bigg[ {b_n^{\nu+1} \over a_{n}^{\nu+1} } -  { b_{n+1}^\nu \over  a_{n+1}^{\nu} } - n \bigg].\eqno(3.24)$$
On the other hand, combining with (3.23), (3.24),  we find

$$ \lambda \left[\gamma^\nu_{n,n} \right]^2 + {1\over \lambda} \bigg[ { A_{n}^{\nu+1} \over \gamma^\nu_{n,n}}\bigg]^2 = B_n^\nu.\eqno(3.25)$$
Moreover, writing (3.23) for $n+1$, we observe one more recurrence relation via (3.24) 

$$\lambda \left[\gamma^\nu_{n,n} \right]^2 +  {1\over \lambda} \bigg[ { A_{n+1}^{\nu+1}  \over \gamma^\nu_{n+1,n+1} }\bigg]^2  = 1 + B_n^{\nu+1} .\eqno(3.26)$$ 
But from (3.22) we have also

$$\lambda \left[\gamma^\nu_{n,n} \right]^2 +  {1\over \lambda} \bigg[ { A_{n+1}^{\nu}  \over \gamma^\nu_{n,n} }\bigg]^2  = 1 + B_n^{\nu+1} .\eqno(3.27)$$ 
Consequently, (3.22), (3.25), (3.26), (3.27) yield the  relations

$$ \left[\gamma^\nu_{n,n-1} \right]^2 + \left[\gamma^\nu_{n,n} \right]^2 = {B_n^\nu \over \lambda} ,\eqno(3.28)$$

$$ \left[\gamma^\nu_{n+1,n} \right]^2 + \left[\gamma^\nu_{n,n} \right]^2 = {1\over \lambda} \left( 1+B_n^{\nu+1} \right),\eqno(3.29)$$

$$ \left[\gamma^\nu_{n,n} \right]^2 = {1\over \lambda}\  \frac {  \left[ A_{n}^{\nu+1}  \right]^2-  \left[ A_{n+1}^{\nu}  \right]^2}{ B_n^\nu - B_n^{\nu+1}  -1},\eqno(3.30)$$

$$ \left[\gamma^\nu_{n,n} \right]^2 = {1\over \lambda} \  \frac{  \left[ A_{n+1}^{\nu}  \right]^2 B_n^\nu - \left[ A_{n}^{\nu+1}  \right]^2 \left[B_n^{\nu+1}  +1\right] } {  \left[ A_{n+1}^{\nu}  \right]^2-  \left[ A_{n}^{\nu+1}  \right]^2},\eqno(3.31)$$

$$\bigg[ \left[ A_{n}^{\nu+1}  \right]^2-  \left[ A_{n+1}^{\nu}  \right]^2\bigg]^2 + \bigg[  \left[ A_{n+1}^{\nu}  \right]^2 B_n^\nu- \left[ A_{n}^{\nu+1}  \right]^2 \left[B_n^{\nu+1}  +1\right] \bigg]$$

$$\times \bigg[  B_n^\nu - B_n^{\nu+1}  -1 \bigg] = 0.\eqno(3.32)$$
Finally in this section we analyse the following series with the backward  difference (2.15)

$$\sum_{k=0}^\infty \nabla \bigg[ P^\nu_n(k; t,\lambda)  \rho_{k+\nu+1}(t) { \lambda^k \over k!} \bigg]  (-k)_j .\eqno(3.33)$$ 
Invoking summation by parts, we represent (3.33) in the form

 $$\sum_{k=0}^\infty \nabla \bigg[ P^\nu_n(k; t,\lambda)  \rho_{k+\nu+1}(t) { \lambda^k \over k!} \bigg]  (-k)_j = -  \sum_{k=0}^\infty P^\nu_n(k; t,\lambda)  \rho_{k+\nu+1}(t) { \lambda^k \over k!}  \Delta  \bigg[  (-k)_j \bigg].$$
Working out the right-hand side of the latter equality by virtue of (1.13), we observe that it is equal to $0$ for $j=0,1,\dots, n-1.$ Therefore

$$ \sum_{k=0}^\infty \nabla \bigg[ P^\nu_n(k; t,\lambda)  \rho_{k+\nu+1}(t) { \lambda^k \over k!} \bigg]  (-k)_j = 0,\quad  j=0,1,\dots, n-1.\eqno(3.34)$$
When $j=n$ we write via (1.13) 

$$\sum_{k=0}^\infty P^\nu_n(k; t,\lambda)  \rho_{k+\nu+1}(t) { \lambda^k \over k!}  \Delta  \bigg[  (-k)_n \bigg] = (-1)^n \sum_{k=0}^\infty P^\nu_n(k; t,\lambda)  \rho_{k+\nu+1}(t)  \lambda^k \bigg[  {k+1\over (k+1-n)!} - {1\over (k-n)!} \bigg]$$

$$=  (-1)^n n \sum_{k=0}^\infty P^\nu_n(k; t,\lambda)  \rho_{k+\nu+1}(t) { \lambda^k \over (k+1-n)!} = 0.$$
Consequently, we established the following orthogonality conditions

$$ \sum_{k=0}^\infty \nabla \bigg[ P^\nu_n(k; t,\lambda)  \rho_{k+\nu+1}(t) { \lambda^k \over k!} \bigg]  (-k)_j = 0,\quad  j=0,1,\dots, n,\eqno(3.35)$$
i.e. due to (1.10)

$$ \sum_{k=0}^\infty \bigg[ \bigg( \lambda (k+\nu) P^\nu_n(k; t,\lambda) - k P^\nu_n(k-1; t,\lambda) \bigg) \rho_{k+\nu}(t) \bigg.$$

$$\bigg. +\lambda t P^\nu_n(k; t,\lambda) \rho_{k+\nu-1} (t) \bigg]   { \lambda^k \over (k-j)!} = 0,\quad  j=0,1,\dots, n.\eqno(3.36)$$

\section{ The modified sequence of polynomials for the weight $t^{-k} \omega_k(\lambda t)$}

This last section is concerned with the sequence of polynomials $\{Q^\nu_n\}_{n\ge 0}$,  being  orthogonal with respect to the weight $\rho_{k+\nu+1}(\lambda t)  \lambda^k/ k! $, i.e.

$$\sum_{k=0}^\infty Q^\nu_n(k; t,\lambda) Q^\nu_m(k; t,\lambda) \rho_{k+\nu+1}(\lambda t) { \lambda^k\over k!} =  \delta_{n,m}.\eqno(4.1)$$ 
The orthogonality conditions can be written accordingly,

$$\sum_{k=0}^\infty Q^\nu_n(k; t,\lambda)  \rho_{k+\nu+1}(\lambda t) { \lambda^k k^j\over k!} = 0,\ j =0,\dots, n-1,\eqno(4.2)$$ 

$$\sum_{k=0}^\infty Q^\nu_n(k; t,\lambda)  \rho_{k+\nu+1}(\lambda t) { \lambda^k (-k)_j  \over k!} = 0,\ j =0,\dots, n-1,\eqno(4.3)$$ 
and explicitly  

$$Q_n^\nu(x; t,\lambda) = \alpha^\nu_n x^n + \beta^\nu_n x^{n-1} + \hbox{lower degrees}.\eqno(4.4)$$
Denoting coefficients of the three term recurrence relation  by $q^\nu_n, h^\nu_n$, it has

$$ x Q_n^\nu(x; t,\lambda) = q^\nu_{n+1} Q_{n+1}^\nu(x; t,\lambda) + h^\nu_n Q_n^\nu(x; t,\lambda) + q^\nu_n Q_{n-1} ^\nu(x; t,\lambda),\eqno(4.5)$$
where
$$ q^\nu_n= {\alpha^\nu_{n-1} \over \alpha^\nu_n},\quad\quad h^\nu_n= {\beta^\nu_{n} \over \alpha^\nu_n} - {\beta^\nu_{n+1} \over \alpha^\nu_{n+1} }.$$
Then we differentiate (4.2) by $t$ which is allowed under the same conditions as in Section 2 to derive via  (1.9) the equalities

$$\sum_{k=0}^\infty  {\partial Q^\nu_n  \over \partial t} (k; t,\lambda)  \rho_{k+\nu+1}(\lambda t) { \lambda^k k^j \over k!}  -  \lambda \sum_{k=0}^\infty  Q^\nu_n  (k; t,\lambda)  \rho_{k+\nu}(\lambda t) { \lambda^k k^j \over k!} = 0,\ j =0,\dots, n-1.\eqno(4.6)$$ 
Meanwhile, the differentiation by $\lambda$ implies

$$ \sum_{k=0}^\infty  {\partial Q^\nu_n  \over \partial \lambda} (k; t,\lambda)  \rho_{k+\nu+1}(\lambda t) { \lambda^k k^j \over k!}  -  t \sum_{k=0}^\infty  Q^\nu_n  (k; t,\lambda)  \rho_{k+\nu}(\lambda t) { \lambda^k k^j \over k!}$$

$$+  {1\over \lambda} \sum_{k=0}^\infty  Q^\nu_n  (k; t,\lambda)  k \rho_{k+\nu+1}(\lambda t) { \lambda^k k^j \over k!}  = 0,\quad  j =0,\dots, n-1.\eqno(4.7)$$ 
The latter term in (4.7) can be rewritten by virtue of the orthogonality and the three term recurrence relation, and we obtain via (4.5) 

$$ \sum_{k=0}^\infty  Q^\nu_n  (k; t,\lambda)  k \rho_{k+\nu+1}(\lambda t) { \lambda^k k^j \over k!} =   q_n \sum_{k=0}^\infty  Q^\nu_{n-1}  (k; t,\lambda)  \rho_{k+\nu+1}(\lambda t) { \lambda^k k^j \over k!}, \quad  j =0,\dots, n-1. $$
Therefore, recalling (4.6), (4.7), we establish the following orthogonality conditions

$$\sum_{k=0}^\infty \bigg[  t {\partial Q^\nu_n  \over \partial t} (k; t,\lambda) - \lambda {\partial Q^\nu_n  \over \partial \lambda} (k; t,\lambda) -  q_n Q^\nu_{n-1}  (k; t,\lambda)\bigg] \rho_{k+\nu+1}(\lambda t) { \lambda^k k^j \over k!}  = 0,\ j =0,\dots, n-1.\eqno(4.8)$$ 
Hence, owing to the uniqueness, we establish the partial differential-difference equation for orthogonal polynomials $Q_n^\nu$

$$  t {\partial Q^\nu_n  \over \partial t} (x; t,\lambda) - \lambda {\partial Q^\nu_n  \over \partial \lambda} (x; t,\lambda) -  q_n Q^\nu_{n-1}  (x; t,\lambda) = {1\over \alpha_n^\nu } \left[ t {\partial \alpha^\nu_n  \over \partial t}  - \lambda {\partial \alpha^\nu_n  \over \partial \lambda} \right] Q^\nu_n (x; t,\lambda).\eqno(4.9)$$
But since similar to (3.4) one finds

$$ 2 \sum_{k=0}^\infty  Q^\nu_n  (k; t,\lambda)  {\partial Q^\nu_n  \over \partial t} (k; t,\lambda)  \rho_{k+\nu+1}(\lambda t) { \lambda^k \over k!} - \lambda   \sum_{k=0}^\infty  \left[Q^\nu_n  (k; t,\lambda) \right]^2  \rho_{k+\nu}(\lambda t) { \lambda^k  \over k!}  = 0,$$

$$ 2 \sum_{k=0}^\infty  Q^\nu_n  (k; t,\lambda)  {\partial Q^\nu_n  \over \partial \lambda} (k; t,\lambda)  \rho_{k+\nu+1}(\lambda t) { \lambda^k \over k!} - t  \sum_{k=0}^\infty  \left[Q^\nu_n  (k; t,\lambda) \right]^2  \rho_{k+\nu}(\lambda t) { \lambda^k  \over k!} $$

$$+ {1\over \lambda} \sum_{k=0}^\infty  \left[Q^\nu_n  (k; t,\lambda) \right]^2  k \rho_{k+\nu+1}(\lambda t) { \lambda^k  \over k!}  = 0.$$
Consequently,   we get from (4.5) and the latter equalities 

$${1\over \alpha_n^\nu } \left[ t {\partial \alpha^\nu_n  \over \partial t}  - \lambda {\partial \alpha^\nu_n  \over \partial \lambda} \right]  =  {h^\nu_n \over 2}.\eqno(4.10)$$
So,  equation (4.9) has the form

$$  t {\partial Q^\nu_n  \over \partial t} (x; t,\lambda) - \lambda {\partial Q^\nu_n  \over \partial \lambda} (x; t,\lambda) -  q^\nu_n Q^\nu_{n-1}  (x; t,\lambda) =  {h^\nu_n\over 2} Q^\nu_n (x; t,\lambda).\eqno(4.11)$$
Then,  equating coefficients in front of $x^{n-1}$ in (4.11) and using (4.10), we derive the Toda-type first order partial differential equation for recurrence coefficients in (4.5)

$$ t {\partial h^\nu_n  \over \partial t}  - \lambda {\partial h^\nu_n  \over \partial \lambda} + [q^\nu_{n+1}]^2 - [q^\nu_n]^2 = 0.\eqno(4.12)$$
To complete the Toda-type system we add an analog of the equation (3.9) which can be verified directly

$$ t {\partial [q^\nu_n]^2  \over \partial t}  - \lambda {\partial [q^\nu_n]^2  \over \partial \lambda} =   [q^\nu_n]^2  \bigg[ h_{n-1}^\nu -  h_{n}^\nu \bigg].\eqno(4.13)$$
Finally,  from the differentiation and  orthogonality we get 

$$0=  \sum_{k=0}^\infty    {\partial Q^\nu_n  \over \partial t} (k; t,\lambda) Q^\nu_{n-1}  (k; t,\lambda) \rho_{k+\nu+1}(\lambda t) { \lambda^k \over k!} - \lambda   \sum_{k=0}^\infty  Q^\nu_n  (k; t,\lambda)  Q^\nu_{n-1}  (k; t,\lambda) \rho_{k+\nu}(\lambda t) { \lambda^k  \over k!},$$

$$ 0= \sum_{k=0}^\infty    {\partial Q^\nu_n  \over \partial \lambda} (k; t,\lambda)  Q^\nu_{n-1}  (k; t,\lambda)\rho_{k+\nu+1}(\lambda t) { \lambda^k \over k!} - t  \sum_{k=0}^\infty  Q^\nu_n  (k; t,\lambda) Q^\nu_{n-1}  (k; t,\lambda)  \rho_{k+\nu}(\lambda t) { \lambda^k  \over k!} $$

$$+ {1\over \lambda} \sum_{k=0}^\infty  Q^\nu_n  (k; t,\lambda) Q^\nu_{n-1}  (k; t,\lambda)  k \rho_{k+\nu+1}(\lambda t) { \lambda^k  \over k!}.$$
This implies immediately the recurrence relation 

$$   { \beta_n^\nu \over \alpha_{n-1}^\nu \alpha_n^\nu} \bigg[ \lambda {\partial \alpha^\nu_n  \over \partial \lambda}  - t {\partial \alpha^\nu_n  \over \partial t}\bigg]  - q_n^\nu = 0,$$
i.e., owing to (4.10),  we end up with the equality (cf. (3.12)) 

$$2 [ q_n^\nu ]^2 + {h^\nu_n \beta^\nu_n \over  \alpha_{n}^\nu} = 0.\eqno(4.14)$$

\bibliographystyle{amsplain}

\end{document}